\documentclass[a4paper 12pt]{article}
\newtheorem{Theorem}{\bf Theorem}
\newtheorem{lemma}{\bf Lemma}
\newtheorem{corollary}{\bf Corollary}
\begin{document}
\author{{\bf Agata Smoktunowicz}}
\date{Institute of Mathematics, Polish Academy of Sciences,
00-956 Warsaw 10, \'{S}niadeckich 8, P.O.Box 21, Poland.}
\title{On primitive ideals in polynomial rings over nil rings}
\maketitle
\begin{abstract} Let $R$ be a nil ring. We prove that primitive ideals
in the polynomial ring $R[x]$ in one indeterminate
over $R$ are of the form $I[x]$ for some ideals $I$ of $R$.
\end{abstract}$ $

$ $
All considered rings are associative but not necessarily have
identities. K{\" o}the's conjecture states
that a ring without nil
ideals has no one-sided nil ideals. It
is equivalent $[4]$ to the assertion that polynomial rings over
nil rings are Jacobson radical.
Our main result states that if $R$ is a nil ring
 and $I$ an ideal in $R[x]$ (the polynomial
 ring in one indeterminate over $R$) then $R[x]/I$
is Jacobson radical if and only if $R/I'[x]$
is Jacobson radical, where $I'$ is the
ideal of $R$ generated by coefficients of
polynomials from $I$.
Also if $R$ is a nil ring and $I$ is a primitive
ideal of $R[x]$ then $I=M[x]$ for some ideal $M$ of $R$.
 It was asked by Beidar, Fong and Puczy{\l}owski $[1]$
 whether polynomial rings over nil rings are not (right ) primitive.
We show that affirmative answer to this question is equivalent to
the K{\" o}the conjecture.
 We also answer in the negative Question 2 from $[1]$ (Corollary 1).

 It is known that
  if a polynomial ring $R[x]$ is primitive
then $R$ need not be primitive $[3]$
( see also Bergman's example in $[5]$).
   Let $R$ be a prime ring and $I$ a nonzero ideal of $R$.
 Then $R$ is a primitive ring if and only if $I$ is a
 primitive ring [6].
Since the Hodges example has a nonzero Jacobson radical
 it follows that  polynomial rings over Jacobson radical
 rings can be right and left primitive (see also
 Theorem 3).

We recall some definitions after $[9]$ (see also $[2]$, $[5]$).

A  right ideal of a ring $R$ is
called {\em modular } in $R$ if
and only if there exists an
element $b\in R$ such that $a-ba\in Q$
for every $a\in R$.
 If $Q$ is a modular maximal right ideal of $R$ then
 for every $r\notin Q$, $rR+Q=R$.

 An ideal $P$ of a ring $R$ is right primitive  in $R$
if and only if there exists a modular maximal right
ideal $Q$ of $R$ such that $P$ is the maximal ideal contained in $Q$.

In this paper $R[x]$ denote the polynomial ring in
one indeterminate over $R$. Given polynomial
 $g\in R[x]$ by $\deg (g)$ we denote the degree of $R$, i.e.,
 the minimal number $d$ such that $g\in R+Rx+\ldots +Rx^{d}$.
Given $a\in R$ let $<a>_{R}$ denote the ideal generated by $a$ in $R$ and
$ <a>_{R[x]}$ the ideal generated by $a$ in $R[x]$. For a ring $R$, by $J(R)$ we
 denote the Jacobson radical of $R$.

 We write $I\triangleleft  R$ if
 $I$ is a two-sided ideal of a
 ring $R$, and $Q\triangleleft _{r}R$ if $Q$ is a right ideal of $R$.
\begin{lemma}
 Let $R$ be a ring , $r\in R$, $f\in R[x]$,
  $Q\triangleleft _{r} R[x]$ and $b\in R[x]$ be such that
  $a-ba\in Q$ for every  $a\in R[x]$.
  Suppose $b -xrf\in Q$. Then for every $i$
  there is $f_{i}\in R[x]$
such that $b-x^{i}rf_{i}\in Q$ and $\deg f_{i}\leq
 \deg f$ for all $i$.
\end{lemma}
{\bf Proof.} We proceed by induction on $n$. If $n=1$ we put
$f_{1}=f$. Suppose Lemma
 holds for some $n\geq 1$
 with $f_{n}=a_{0}+\ldots +a_{k}x^{k}$ ($k\leq \deg f)$.
Then Lemma holds for $f_{n+1}=fra_{0}
+\sum_{i=1}^{k}a_{i}x^{i-1}$.
\begin{lemma}
Let $R$ be a ring and $I\triangleleft R[x]$
 and $a_{0}+a_{1}x+\ldots +a_{k}x^{k}\in I$,
$a_{0},\ldots , a_{k}\in R$. Denote $U=<a_{k}>_{R[x]}$.
\begin{itemize}
\item[1.] Suppose $h\in  {U}^{l}$ for some $l\geq 1$
 and $\deg h\geq k$.
 Then there is $g\in {U}^{l-1}$ such that ,
$h-g\in I$ and $\deg g<\deg h$.
\item[2.]
Suppose $Q\triangleleft _{R} R[x]$, $I\subseteq Q$.
Let $b\in R[x]$ be such that $a-ba\in Q$ for every
$a\in R[x]$. If $b-xrf\in Q$ and $b-g\in Q$ where $f, g\in R[x],
g\in U^{\deg f}$ then for every $i>\deg g$ there is $g_{i}\in R[x]$
such that $b-x^{i}rg_{i}\in Q$ and $\deg g_{i}<k$.
\end{itemize}
\end{lemma}
{\bf Proof.} {\bf 1.} Let $h=\sum_{i=0}^{t}c_{i}x^{i}$,
$t=\deg h$. Since $h\in U^{l}$ then
$c_{t}=\sum_{i}p_{i}a_{k}q_{i}$ for some $p_{i},q_{i}\in R\cup \{1\}$,
 $q_{i}\in U^{l-1}$.
 We put $g=h-c_{t}x^{t}-
\sum_{i}p_{i}(\sum_{j=0}^{k-1}a_{j}x^{j})q_{i}x^{t-k}$.\\
{\bf 2. }Let $f_{i}$ be as in Lemma 1 applied
for $f$. Let $g=\sum_{i=0}^{t}c_{i}x^{i}$, $c_{i}\in R$,
$c_{i}\in U^{\deg f}$.
 For a natural $n>t$, we put
 $h_{n}=\sum_{i=0}^{t}f_{n-i}c_{i}$. Observe that
 $\deg h_{n}\leq \deg f$, $h_{n}\in U^{\deg f}$
 and $b-x^{n}rh_{n}\in Q$.
 Applying several times Lemma 2.1
 for $h=h_{n}$, we get that there are
$g_{n}\in R[x]$ such that $b-x^{n}rg_{n}\in Q$
for all $n>t$ and $\deg g_{n}\leq k$.
We are done.

Let $R$ be a ring.
Given $r\in R$ and $Q\triangleleft _{r} R$
 and $b\in R[x]$ such that $a-ba\in Q$ for every
 $a\in R[x]$
 we say that $v$ is a ``good number for $r$`` if
for all sufficiently
large $n$, there are $f_{n}\in R[x]$
such that $\deg f_{n}\leq v$ and $b-x^{n}rf_{n}\in Q$.
 \begin{lemma}
Let $R$ be a ring, $Q\triangleleft_{R} R[x]$,
 and $b\in R[x]$ be such that $a-ba\in Q$ for every $a\in R[x]$.
 Let $r\in R$, $r\notin Q$
and let $v$ be a good number for all $a\in rR$, $a\notin Q$.
Suppose there are $p,p'\notin Q$,
$p,p'\in rR$  such that
 $(pR+Q)\cap (p'R+Q)\subseteq Q$. Then $v-1$ is
 a good number for $r$.
\end{lemma}
{\bf Proof. }Since $v$ is a good number for $p,p'$
 then for sufficiently large $n$ there are $g_{n}\in pR[x]$,
$g_{n}'\in p'R[x]$ with $\deg g_{n}, \deg g_{n}'\leq v$
 such that $b-x^{n}g_{n}\in Q$ and $b-x^{n}g_{n}'\in Q$.
Let $g_{n}=p_{n,0}+p_{n,1}x+\ldots +p_{n,v}x^{v}$,
$g_{n}'=p_{n,0}'+p_{n,1}'x+\ldots +p_{n,v}'x^{v}$,
 where all $p_{n,i}\in pR, p_{n,i}'\in p'R$.
If for all sufficiently large $n$ either $p_{n,v}\in Q$
or $p_{n,v}'\in Q$ then $v-1$ is a good number for $r$.
Suppose that there is $m$ such that $p_{m,v}\notin Q$
 and $p_{m,v}'\notin Q$.
 Since $p_{m,v}\in pR$, $p_{m,v}'\in p'R$
 then $p_{m,v}-p_{m,v}'\notin Q$ since
 otherwise $p_{m,v}\in (pR+Q)\cap (p'R+Q)$
 contrary our assumptions.
Since $v$ is a good number for $c=p_{m,v}-p_{m,v}'$
then for sufficiently large $n$ there
are $h_{n}\in R[x]$,
$\deg h_{n}\leq v$, $b-x^{n}ch_{n}\in Q$. Now
$ch_{n}=c\bar {g}_{n}+cr_{n}x^{v}$ for some
$r_{n}\in R$, $\bar {g}_{n}\in R[x]$,
$\deg \bar {g}_{n}\leq v-1$. Thus (for
sufficiently large $n$) $b-x^{n}k_{n}\in Q$
 where $k_{n}=ch_{n}+(g_{m}'-g_{m})r_{n} =
 c\bar {g}_{n}+\sum_{i=0}^{v-1}(p_{m,i}'-p_{m,i})x^{i}r_{n}$.
 Note that $k_{n}\in rR[x]$ since $p,p'\in rR$. Since
 $\deg k_{n}\leq v-1$ then $v-1$
 is a good number for $r$.
\begin{Theorem}
Let $R$ be a nil ring and $I$ be a primitive  ideal
 in $R[x]$ (if any). Then $I=I'[x]$ for some ideal $I'$ in $R$.
\end{Theorem}
{\bf Proof. }Suppose the contrary, that there are $a_{0}, a_{1},\ldots ,
a_{k}\in R$, $a_{k}\notin I$ such that $a_{0}+a_{1}x+\ldots +a_{k}x^{k}\in I$.
 From the definition there
 is a modular maximal
 right ideal $Q$ in $R[x]$ such that
$I$ is the maximal possible ideal contained in $Q$.
Now for every $r\notin Q$, $r\in R[x]$, $rR[x]+ Q=R[x]$.
Since $Q$ is modular, there is $b\in R[x]$
 such that $a-ba\in Q$ for every $a\in R[x]$.
  Observe first that if $r\notin Q$ then $rx\notin Q$ (otherwise
  $xR\subseteq Q, xR\subseteq I$, impossible since $R$ is nil). Thus
  there is $f\in R[x]$ such that
  $b-xrf\in Q$. Denote $U=<a_{k}>_{R[x]}$. Since $I$ is a primitive ideal
  it is prime. Consequently $U^{\deg f}\not\subseteq I$,
  $U^{\deg f}\not\subseteq Q$.
 Hence $b-g\in Q$ for some $g\in U^{\deg f}$.
By Lemma 2.2 we get that $k-1$ is a good number for all
$r\in R$, $r\notin Q$. Let $v$ be minimal such that
$v$ is a good number for all $r\in R$,
$r\notin Q$. We will show that
 $v=0$ and hence get a contradiction (since $R$ is a nil ring ).
 Suppose $v>0$. Let $r\in R$, $r\notin Q$. We will show that $v-1$
is a good number for $r$. Since $v$ is a good number for
$r$ then for some $i$ there are $f_{i}, f_{i+1}, \ldots , f_{i+k}\in rR[x]$
 such that $b-x^{j}f_{j}\in Q$ where $i\leq j\leq i+k$.
 Now each $f_{j}$ can be written as $f_{j}=g_{j}+x^{v}c_{j}$,
$c_{j}\in rR$, $g_{j}\in R[x]$,
 $\deg g_{j}<v$. Let $e_{j}=c_{j}^{n_{j}}$ where $n_{j}$
 is minimal possible such that
$c_{j}^{n_{j}}\notin Q$ where $i\leq j\leq i+k$.
(if $c_{j}\in Q$ we put $e_{j}=1$).
 Now either $v-1$ is a good number for $r$  or
 there is  $s\notin Q$ such that
 $s\in \bigcap _{j=i}^{i+k}(e_{j}R+Q)$
  and $s\in e_{i}R$ (by Lemma 3).
 Then
 $s-e_{j}d_{j}\in Q$ for some $d_{j}\in R$, $i\leq j\leq i+k$.
 Since $v$ is a good number for $s$ then
 for sufficiently large
 $n$ there are $\bar {f}_{n}\in sR[x]$,
$\deg \bar {f}_{n}\leq v $ such that $b-x^{n}\bar {f}_{n}\in Q$.
Let $\bar {f}_{n}=\sum_{j=0}^{v} sb_{j}x^{j} $. Thus $\bar
{f}_{n}-\sum_{j=0}^{v} e_{i+v-j}d_{i+v-j}b_{j}x^{j} \in Q$. Now
since  $b-x^{j}f_{j}\in Q$ then
 $(b-x^{j}f_{j})e_{j}d_{j}\in Q$.
 Thus $\bar {f}_{n}-\bar {g}_{n}\in Q$ where
 $\bar {g}_{n}=\sum_{j=0}^{v} x^{i+v-j}f_{i+v-j}e_{i+v-j}d_{i+v-j}b_{j}x^{j}$.
 Note that $\bar {g}_{n}=
x^{i+v}\sum_{j=0}^{v}f_{i+v-j}e_{i+v-j}d_{i+v-j}b_{j}$.
 Observe that there are $t_{l}\in rR[x]$, $\deg t_{l}\leq v-1$
  such that $t_{l}-f_{l}e_{l}\in Q$
 for $i\leq l\leq i+v$.
 Denote $h_{n}=\sum_{j=0}^{v}t_{i+v-j}d_{i+v-j}b_{j}$.
 Then $h_{n}\in rR[x]$ and $\deg h_{n}\leq v-1$.
 Since $b-x^{n}\bar {f}_{n}\in Q$ then
 $b-x^{n}\bar {g}_{n}\in Q$. Thus $b-x^{i+v+n}h_{n}\in Q$.
 Since it holds for all sufficiently large $n$
 we get that $v-1$ is a good number for $r$.
\begin{Theorem}
Let $R$ be a nil ring and $I\triangleleft R[x]$ and let  $\bar
I\triangleleft R$ be the ideal of $R$ generated by coefficients of
polynomials from $I$. Then $R[x]/I$ is Jacobson radical if and
only if $R[x]/\bar {I}[x]$ is Jacobson radical.
\end{Theorem}
{\bf Proof. }Suppose the contrary, that $R/\bar {I}$ is Jacobson
radical and $R[x]/I$ is not Jacobson radical. Then there is a
primitive ideal $P$ of $R[x]/I$ such that $P\neq R[x]/I$. Now
$P+I$ is a primitive ideal in $R[x]$. Thus $P+I=\bar {P}[x]$ for
some ideal $\bar {P}$ in $R$ by Theorem 1. Since $P+I\neq R[x]$
and $\bar {I}\subseteq \bar {P}$ then $R[x]/\bar {I}[x]$ is not
Jacobson radical, a contradiction. The other inclusion is clear.
\begin{corollary}
 If $N$ is a nil ring then the polynomial ring $N[x]$ can not be
 homomorphically mapped onto a simple primitive ring.
\end{corollary}
 Krempa [4] showed that the K{\" o}the conjecture is equivalent to the
assertion that polynomial rings over nil rings
are Jacobson radical. From this and Theorem 1 we get:
\begin{corollary}
 The K{\" o}the conjecture is equivalent to the statement
 ``polynomial rings (in one indeterminate) over nil rings are
 not right primitive``.
\end{corollary}

  Simple Jacobson radical but not nil rings were constructed in [7, 8].
  (Rings in [7, 8] are not nil since they satisfy the relation
  x=yxxy).

\begin{Theorem}
 Let $R$ be a simple Jacobson radical ring which is not nil. Then
the polynomial ring $R[x]$ in one indeterminate
over $R$ is right primitive.
\end{Theorem}
{\bf Proof. }Since $R$ is not nil there is  $b$ in $R$
 such that $b$ is not nilpotent.
    Let $Q$ be a right ideal in $R[x]$ maximal with the property that
    $xb\notin Q$ and $r-xbr\subseteq Q$ for every $r\in R[x]$.
    Then $Q$ is a maximal modular right
    ideal in $R[x]$. We will show that if $I$
   is a two-sided ideal of $R[x]$ and $I\subseteq Q$
    then $I=0$. Let $a_{0}+a_{1}x+\ldots a_{n}x^{n}\in I$ where
     $a_{0},\ldots ,a_{n}\in R, a_{n}\neq 0$.
 Since $R$ is a simple ring $n>0$ and there are $b_{j},c_{j}\in R$,
  $j=1,2,\ldots ,m$
  such that $\sum_{j=1}^{m}b_{j}a_{n}c_{j}=b$. Denote
$g=d_{0}+\ldots +d_{n}x^{n}$, where
$d_{i}=\sum_{j=1}^{m}b_{j}a_{i}c_{j}$ for $0\leq i\leq n$.
 Let $d=\sum_{k=1}^{n-1}b^{n-1-k}d_{k}$.
Note that $r-(bx)^{k}r\in Q$ for every $r\in R[x]$, $k>0$.
 Thus $d-\sum_{k=1}^{n}b^{n-1}x^{k}d_{k}\in Q$.
 Since $g\in I$ then $b^{n-1}g\in Q$. Consequently
 $d +b^{n}x^{n}\in Q$ (since $d_{n}=b$).
 Thus $dr+r\in Q$ for every $r\in R[x]$. Since $d\in R$ and $R$ is Jacobson radical and $d\in R$ then $d+dr+r=0$ for some $r\in R$.
  Thus $d\in Q$, impossible.\\

 \noindent
{\Large\bf REFERENCES:}
\begin{itemize}
\item[1]K.Beidar, Y.Fong, E.R.Puczy{\l}owski,
Polynomial rings over nil rings cannot be homomorphically
mapped onto rings with nonzero idempotents, J.Algebra, 238 (2001) 389-399.
\item[2]N.J.Divinsky, Rings and Radicals, Allen, London, 1965.
\item[3]T.J.Hodges, An example of a primitive polynomial ring,
J.Algebra 90 (1984) 217-219.
\item[4]J.Krempa, Logical connections among some open
problems concerning nil rings, Fund. Math. 76 (1972) 121-130.
\item[5]T.Y.Lam, Exercises in Classical Ring Theory,
Problem Books in Mathematics, Berlin-Heidelberg-New
York, Springer-Verlag, 1995 pp. 122-123.
\item[6]C.Lanski, R.Resco, L.Small, On the primitivity of
prime rings, J.Algebra 59 (1979) 395-398.
\item[7]E.Sasiada, Solution of the problem  of
existence of  a simple radical ring,
Bull.Acad.Polon.Sci.Ser.Sci.Math.Astr.Phys., 9 (1961) 257.
\item[8]E.Sasiada and P.M.Cohn, An example of a simple
radical ring, J.Algebra, 5 (1967) 373-377.
\item[9]F.A.Sz{\' a}sz, Radicals in rings, Akademia Kiado,
Budapest, 1981.
\end{itemize}
\end{document}